\documentclass[12pt, 14paper,reqno]{amsart}
\usepackage{amsmath,amsfonts,amssymb}
\usepackage[breaklinks]{hyperref}
\usepackage{graphicx}
\usepackage{longtable}
\makeatletter
\@namedef{subjclassname@2020}{%
  \textup{2020} Mathematics Subject Classification}
\makeatother

        \makeatletter
\def\legendre@dash#1#2{\hb@xt@#1{%
  \kern-#2\p@
  \cleaders\hbox{\kern.5\p@
    \vrule\@height.2\p@\@depth.2\p@\@width\p@
    \kern.5\p@}\hfil
  \kern-#2\p@
  }}
\def\@legendre#1#2#3#4#5{\mathopen{}\left(
  \sbox\z@{$\genfrac{}{}{0pt}{#1}{#3#4}{#3#5}$}%
  \dimen@=\wd\z@
  \kern-\p@\vcenter{\box0}\kern-\dimen@\vcenter{\legendre@dash\dimen@{#2}}\kern-\p@
  \right)\mathclose{}}
\newcommand\legendre[2]{\mathchoice
  {\@legendre{0}{1}{}{#1}{#2}}
  {\@legendre{1}{.5}{\vphantom{1}}{#1}{#2}}
  {\@legendre{2}{0}{\vphantom{1}}{#1}{#2}}
  {\@legendre{3}{0}{\vphantom{1}}{#1}{#2}}
}
\def\dlegendre{\@legendre{0}{1}{}}
\def\tlegendre{\@legendre{1}{0.5}{\vphantom{1}}}
\makeatother
\theoremstyle{plain}
\theoremstyle{plain}

\newtheorem{rem}{Remark}[section]
\numberwithin{equation}{section}
\newtheorem{thm}{Theorem}[section]
\newtheorem{lem}{Lemma}[section]
\newtheorem{cor}{Corollary}[section]
\newtheorem{prop}{Proposition}[section]

\numberwithin{equation}{section}

\usepackage[breaklinks]{hyperref}
\makeatletter
\let\oldHyPsd@CatcodeWarning\HyPsd@CatcodeWarning
\renewcommand{\HyPsd@CatcodeWarning}[1]{
  \ifnum\pdfstrcmp{#1}{math shift}=0    
  \else                                 
    \oldHyPsd@CatcodeWarning{#1}
  \fi
}
\makeatother

\begin{document} 
\title[Elliptic Curve]{on the family of elliptic curves $y^2=x^3-5pqx$}
\author{Arkabrata Ghosh}
\address{G. Arkabrata@ Ashoka University, National Capital Region P.O, Plot no 2, Rajiv Gandhi Education City, Rai, Sonipat, Haryana 131029}
\email{arka2686@gmail.com}



\maketitle

\section*{Abstract}

This article considers the family of elliptic curves given by $E_{pq}: ^2=x^3-5pqx$ and certain conditions on odd primed $p$ and $q$. More specifically, we have proved that if $p \equiv 33 \pmod {40}$ and $ q \equiv 7 \pmod {40}$, then the rank of $E_{pq}$  is zero over both $ \mathbb{Q} $ and $ \mathbb{Q}(i) $. Furthermore, if the primes $ p $  and $q$ are of the form $ 40k + 33 $ and $ 40l + 27$, where $k,l \in \mathbb{Z}$ such that $(25k+ 5l +21)$ is a perfect square, then the given family of elliptic curves has rank one over $\mathbb{Q}$ and rank two over $\mathbb{Q}(i)$. Finally, we have shown that torsion of $E_{pq}$ over $\mathbb{Q}$ is isomorphic to $ \mathbb{Z}/ 2\mathbb{Z}$.

\section{Introduction}
The arithmetic of the elliptic curve is one of the most fascinating branches of mathematics as it connects number theory to algebraic geometry. 
Recently, there has been a surge among mathematicians in studying the families of the elliptic curve given by $E: y^2= x^3 + bx $ where $ b \in \mathbb{Z}$ and several people are now working in this direction. In 2007, Spearman \cite{SP07} calculated the values of the prime $p$ for which the elliptic curve $ E: y^2 =x^3 -px$ has rank two.  In the same year, he \cite{Sp07} also discovered the condition on $2p$ for which the elliptic curve $y^2= x^3- 2px$ has rank three.  tions, the rank of the given elliptic curve could be two, three or four. Fujita and Terai \cite{FT11}, in 2011,  considered the elliptic curve of the form $y^2 = x^3- p^{k}x $. where $p$ is prime and $ k = 1, 2, 3$ and found a necessary and sufficient condition for the rank of the given elliptic curve to be equal to one or two. Daghieh and Didari \cite{DD14} determined the rank of the elliptic curves of the form $ y^2 = x^3 -3px$, $ p $ be a prime number. They also studied \cite{DD15} elliptic curves of $y^2 = x^3 -pqx$, where $p$ and $q$ are distinct odd primes. Around the same time, Kim \cite{K15} studied the elliptic curve $y^2 =x^3 \pm 4px$, where $p$ is a prime number. Recently, Mina and Bacani \cite{MB23} studied the curve $y^2= x^3 -3pqx$ for distinct odd primes $p$ and $q$ and find two different sufficient conditions on $p$ and $q$ such that the given elliptic curve, under those conditions, has rank zero and one respectively.

In this article, we aim to find the rank of the following family of elliptic curves given by

\begin{equation}
\label{eq:1.1}
    E_{pq}~ y^2= x^3-5pqx,
\end{equation}

where $p$ and $q$ are distinct odd primes satisfying certain conditions. At first, we will provide conditions on $p$ and $q$ such that the rank of the family elliptic curve $E_{pq}$ given by equation \eqref{eq:1.1} is zero. More precisely, we will prove the following Theorem.

\begin{thm}
    \label{Thm:1.1}
    Suppose $p$ and $q$ be distinct primes satisfying the congruence $ p \equiv 33 \pmod {40}$ and $ q \equiv 7 \pmod {40}$. Then rank of the elliptic curve $E_{pq}: y^2 =x^3-5pqx$ is zero over $\mathbb{Q}$.
\end{thm}

Now we will update the conditions of $p$ and $q$ such that the rank of the elliptic curve given by \eqref{eq:1.1} is exactly one. This is stated in the following Theorem.

\begin{thm}
    \label{thm:1.2} 
   Let $ p = 40k + 33$ and $ q =40 l + 27$ for some $ k, l \in \mathbb{Z}$. If $ (25k + 5l + 21)$ is a perfect square, then the elliptic curve $E_{pq}$ given by \eqref{eq:1.1} has rank one over $\mathbb{Q}$.
\end{thm}

We will end this article by discussing the torsion part of the elliptic curve given by \eqref{eq:1.1}. More precisely, we will prove the following result

\begin{thm}
    \label{thm:1.3}
    Let $p$ and and $q$ be distinct primes satisfying the congruence $ p \equiv 33 \pmod {40}$ and $ q \equiv 7 \pmod {40}$. Then the Torsion subgroup of the elliptic curve $E_{pq}: y^2 =x^3-5pqx$ is isomorphic to $\mathbb{Z}/2\mathbb{Z}$.
\end{thm}

\section{preliminaries}

Let $E$ be an elliptic curve over the field $\mathbb{K}$ of characteristic not equal to $2$ or $3$ and let $ E(\mathbb{K})$ denote the $\mathbb{K}$-rationals points of $E$ over $\mathbb{K}$. Mordell-Weil Theorem asserts that $  E(\mathbb{K})$ is a finitely generated abelian group and can be represented as

\begin{equation*}
    E(\mathbb{K}) \cong \mathbb{Z}^r \oplus E(\mathbb{K})_{tors}, 
\end{equation*}
where $ r \geq 0$ is called the rank of the elliptic curve $E$ over $\mathbb{K}$ and $ E(\mathbb{K})_{tors}$ denotes the torsion sub-group of $ E(\mathbb{K})$ which is a finite abelian group consisting of elements of finite order. The following Theorem \cite{M77} gives an explicit idea about the torsion subgroup $ E(\mathbb{K})_{tors}$ where $\mathbb{K}=\mathbb{Q}$.

\begin{thm}(Mazur's Theorem)
\label{thm:2.5}
    Let $E$ be an elliptic curve defined over  $\mathbb{Q}$. Then 

$$
 E_{tors}(\mathbb{Q}) \cong 
 \begin{cases}  
\mathbb{Z}/ n\mathbb{Z}~~~ ~~~{\ \text{ for }\ 1 \leq n\leq 10, ~\text{or}~ n=12},
\\ 
\mathbb{Z}/ 2 \mathbb{Z} \times \mathbb{Z}/ 2 n \mathbb{Z} {\ \text{ for }\ ~1 \leq n \leq 4}.
\end{cases}
$$ 
\end{thm}

Now to determine the torsion subgroup of the families given by the equation \eqref{eq:1.1},
we need the following Theorem (\cite{ST15}, Page $56$).

\begin{thm}(Nagel-Lutz Theorem)
    \label{thm:2.1}
    We consider the elliptic curve $Y^2= X^3 + aX + b$ with $a, b \in \mathbb{Z}$ and, let $(X,Y) \in E(\mathbb{Q}) $ be a torsion point. Then, 
    \begin{equation*}
       \begin{cases}
         X,Y \in \mathbb{Z}, \\
         \text{Either}~ Y =0, ~\text{or}~ Y^2|(4a^3 + 27b^2). 
        \end{cases}
    \end{equation*}
\end{thm}

Now to compute the rank of the curve given by \eqref{eq:1.1}, we need to use the method of $2$-descent. We will describe it briefly here and one can look at (see \cite{ST15}) for more details. Suppose that $E : Y^2= X^3 + a X^2 + bX$ is an elliptic curve over $\mathbb{Q}$ and $\overline{E}: Y^2 = X^3 - 2aX^2 + (a^2- 4b) X$ is the corresponding isogenous curve to $E$. Hence, there exists an isogeny $ \phi: E \rightarrow \overline{E}$ of degree $2$ given by
\begin{equation*}
    \phi(x,y) ~ = ~ \bigg(\frac{y^2}{x^2}, \frac{y(x^2-b)}{x^2}\bigg).
\end{equation*}
Moreover, let $\mathbb{Q}^{*}$ be the multiplicative group of all non-zero rational numbers, and $\mathbb{Q^*}^{2}$ be its subgroup of squares of elements of $ \mathbb{Q^*}$. Hence, $ \mathbb{Q^*}/\mathbb{Q^*}^{2}$ is the multiplicative group of all non-zero rational numbers modulo squares. Also, we denote set of rational points on $E$ and $\overline{E}$ by $\Gamma$ an $ \overline{\Gamma}$ respectively. Now we define the $2$-descent homomorphism $\alpha: \Gamma \rightarrow \mathbb{Q^*}/\mathbb{Q^*}^{2}$ by

\begin{equation*}
    \alpha(P)= 
               \begin{cases}
                 1 \pmod {\mathbb{Q^*}^2}, ~\text{if}~ P= \mathcal{O}, ~\text{the point of infinity},\\
                 b \pmod {\mathbb{Q^*}^2}, ~\text{if}~ P= (0,0),\\
                  x \pmod {\mathbb{Q^*}^2}, ~\text{if}~ P= (x,y ) ~\text{with}~ x \neq 0.
                 \end{cases}
\end{equation*}

Similarly, we can define the $2$-decent homomorphism on the isogeneous curve $\overline{E}(\mathbb{Q}) $ curve as follows: $\overline{\alpha}:  \overline{\Gamma} \rightarrow  \mathbb{Q^*}/\mathbb{Q^*}^{2}$ by 

\begin{equation*}
    \overline{\alpha}(\overline{P})= 
               \begin{cases}
                 1 \pmod {\mathbb{Q^*}^2}, ~\text{if}~ \overline{P}= \mathcal{\overline{O}}, ~\text{the point of infinity},\\
                 \overline{b} \pmod {\mathbb{Q^*}^2}, ~\text{if}~ \overline{P}= (0,0),\\
                 x \pmod {\mathbb{Q^*}^2}, ~\text{if}~ \overline{P}= (x,y ) ~\text{with}~ x \neq 0,
                 \end{cases}
\end{equation*}

where $\overline{b}= a^2 - 4b$.

The group $\alpha(\Gamma)$ consists of $1,b$ and all factors $b_1$ of $b$, all modulo $\mathbb{Q^*}^2$. Here $b_1 \neq 1, ~\text{or}~ b \pmod {\mathbb{Q^*}^2}$, such that the triple $(N, M, e) \in \mathbb{Z}^3$, where $ M \neq 0, e \neq 0$ solves the Diophantine equation ( or 'torsors')

\begin{equation*}
    \mathcal{T}: N^2 = b_1 M^4 + a M^2 e^2 + b_2 e^4, ~\text{with}~  b_1 b_2 =b,
\end{equation*}

and satisfies the following Criterion :

\begin{equation*}
    gcd(N,e)=gcd(M,e)=gcd(b_1, e) =gcd(b_2, M)= gcd(M,N)=1.
\end{equation*}

Similarly, The group $\alpha(\overline{\Gamma})$ consists of $1, a^2-4b$ and all factors $b_1$ of $a^2-4b$, all modulo $\mathbb{Q^*}^2$. Here $b_1 \neq 1, ~\text{or}~ a^2-4b \pmod {\mathbb{Q^*}^2}$, such that the triple $(N, M, e) \in \mathbb{Z}^3$, where $ M \neq 0, e \neq 0$ solves the Diophantine equation ( or 'torsors')

\begin{equation*}
    \mathcal{T'}: N^2 = b_1 M^4 -2a M^2 e^2 + b_2 e^4, ~\text{with}~  b_1 b_2 =a^2-4b,
\end{equation*}

and the same GCD criterion mentioned above. Now to compute the rank of an elliptic curve $E$, we use the following Proposition(see \cite{ST15}).
\begin{prop}
    \label{prop:2.1} Let $r$ be the rank of $E(\mathbb{Q}) $ and $\alpha$ and $\Bar{\alpha}$ are as above. Then,
    $$
       \frac{1}{4} |\alpha(\Gamma)| \cdot |\overline{\alpha}(\overline{\Gamma})|= 2^r.
    $$
\end{prop}

Now if $\mathbb{K} = \mathbb{Q}(\sqrt{m})$, where $m$ is a square-free integer, we can find the rank of any elliptic curve $E$ over $\mathbb{K}$ by adding ranks of $E$ and its' $m$-twist $E[m]$ over $\mathbb{Q}$. This can be seen from the following result (see \cite{S09}).

\begin{thm}
    \label{thm:2.2}
    Let $\mathbb{K} = \mathbb{Q}(\sqrt{m})$ be a quadratic field, where $m$ is a square-free integer. Let $E: y^2 = x^3 + ax^2 + bx $ be an elliptic curve over $\mathbb{Q}$ and $E[m]: y^2 = x^3 + max^2 + m^2 bx$ be the $m$-twist of $E$. Then

    $$
    rank(E(\mathbb{K}))= rank(E(\mathbb{Q})) + rank(E[m](\mathbb{Q}) ).
    $$
    
\end{thm}
  
\section{Rank of $\texorpdfstring{E_{pq}}{}$}

We will start this section with the following remark.

\begin{rem}
\label{rem:3.1}
Let $p$ and $q$ be two distinct primes. Then the following holds.  \\
(a) If $ p \equiv 33 \pmod {40}$, then $ \genfrac(){}{0}{p}{5}=-1$. \\
(b) If $q \equiv 7 \pmod {40} $, then $ \genfrac(){}{0}{q}{5}=-1$.
\end{rem}

\begin{proof}
    (a) From the definition of the Legendre symbol, we know
    \begin{equation*}
        \genfrac(){}{0}{p}{5} \genfrac(){}{0}{5}{p} =(-1)^{\frac{p-1}{2} \frac{5-1}{2}}= (-1) ^{(p-1)}=1.
    \end{equation*}
    As $ p \equiv 33 \pmod {40} \equiv 3 \pmod 5$, so $ \genfrac(){}{0}{5}{p}=-1$. Hence, from the above equation, we can say that $ \genfrac(){}{0}{p}{5}=-1$.\\
    
    (b) Proceeding similarly as $(a)$, we have
    \begin{equation*}
        \genfrac(){}{0}{q}{5} \genfrac(){}{0}{5}{q} =(-1)^{\frac{p-1}{2} \frac{5-1}{2}}= (-1) ^{(q-1)}=1.
    \end{equation*}
 Now $ q \equiv 7 \pmod{40} \equiv 2 \pmod{5}$, we know $ \genfrac(){}{0}{5}{q}=-1$.  Hence,   $\genfrac(){}{0}{q}{5}=-1$.
\end{proof}

From the equation \eqref{eq:1.1}, we have $  E_{pq}: ~ y^2= x^3-5pqx $ and the corresponding $\overline{E}_{pq}: y^2 = x^3 + 20pqx$. To compute the rank of $E_{pq}$, we first need to determine $|\alpha(\Gamma)|$. Note that $\{ 1, -5pq \} \in \alpha(\Gamma)$ by the definition of $\alpha$. In that case, the following set gives all possible advisors $b_1$ of $-5pq$ modulo $\mathbb{Q^*}^2$.

\begin{equation*}
    S= \{-1, \pm 5q, \pm p, \pm 5p, \pm q, \pm 5, \pm pq, 5pq \}.
\end{equation*}

We then consider the solvability of the following torsors over the set of integers.

\begin{equation*}
  \begin{split}
    \mathcal{T}_1: & N^2 = 5pqM^4 -e^4 \\
    \mathcal{T}_2 : & N^2 = 5qM^4- pe^4 \\
    \mathcal{T}_3 : & N^2 = pM^4- 5qe^4 \\
    \mathcal{T}_4 : & N^2 = pqM^4- 5e^4 \\
    \mathcal{T}_5 : & N^2 = 5M^4- pqe^4 \\
    \mathcal{T}_6 : & N^2 = 5pM^4- qe^4 \\
    \mathcal{T}_7 : & N^2 = qM^4- 5pe^4 \\
    \end{split}
\end{equation*}

\begin{lem}
    \label{le:3.1}
    There are no integer solutions for the torsor $ \mathcal{T}_1: N^2=5pqM^4-e^4$.
\end{lem}

\begin{proof}
    Reducing $ \mathcal{T}_1 $ modulo $q$, we get $N^2 \equiv -e^4 \pmod{q}$. This implies that

    \begin{equation*}
        1 = \genfrac(){}{0}{-e^4}{q}= \genfrac(){}{0}{-1}{q} \genfrac(){}{0}{e^4}{q}= \genfrac(){}{0}{-1}{q}
    \end{equation*}

Now $ \genfrac(){}{0}{-1}{q}=1 \iff q \equiv 1 \pmod 4$. As $ q \equiv 7 \pmod {40} \equiv 3 \pmod 4$, we arrived at a contradiction. Hence we can say that $  \mathcal{T}_1  $ has no solution in $ \mathbb{Z}$.
    
\end{proof}

\begin{lem}
\label{le:3.2}
There are no integer solution for the torsors $ \mathcal{T}_2: N^2=5qM^4-pe^4$ and $ \mathcal{T}_3: N^2=pM^4-5qe^4$.
\end{lem}

\begin{proof}
    
    Reducing $ \mathcal{T}_2$ modulo $5$, we get $N^2 \equiv -pe^4 \pmod 5$. So,

    \begin{equation*}
        1 = \genfrac(){}{0}{-pe^4}{5}= \genfrac(){}{0}{-1}{5} \genfrac(){}{0}{p}{5} \genfrac(){}{0}{e^4}{5} =  \genfrac(){}{0}{p}{5}.
    \end{equation*}
    This is a contradiction as from the remark \eqref{rem:3.1}, 
we know $\genfrac(){}{0}{p}{5}=-1$. Now considering $ \mathcal{T}_3$ modulo $5$  and proceeding similarly as above, we get $ \genfrac(){}{0}{p}{5}=1$, a contradiction.
\end{proof}

\begin{lem}
\label{le:3.3}
There are no integer solution for the torsors $ \mathcal{T}_4: N^2=pqM^4-5e^4$ and $ \mathcal{T}_5: N^2 = 5M^4 -pqe^4$.
\end{lem}

\begin{proof}

As $ p \equiv 33 \pmod {40} \equiv 1 \pmod 4$, we get $\genfrac(){}{0}{-1}{p}=1$. Considering the torsor modulo $p$, we have $ N^2 \equiv -5e^4 \pmod p$. Hence we get, 

\begin{equation*}
    1 = \genfrac(){}{0}{-5e^4}{p}= \genfrac(){}{0}{-1}{p} \genfrac(){}{0}{e^4}{p} \genfrac(){}{0}{5}{p}= \genfrac(){}{0}{5}{p}.
\end{equation*}
    Now $ \genfrac(){}{0}{5}{p}=1 \iff p \equiv 1, 4 \pmod{5}$. As $ p \equiv 33 \pmod {40} \equiv 3 \pmod{5}$, we arrived at a contradiction. Similarly, by reducing $\mathcal{T}_5$ modulo $5$ and using Legendre symbols, we get $  \genfrac(){}{0}{5}{p}=1 $, a contradiction. 
    \end{proof}

\begin{lem}
\label{le:3.4} There are no integer solution for the torsors $ \mathcal{T}_6: N^2=5pM^4-qe^4$ and $ \mathcal{T}_7: N^2 = qM^4 -5pe^4$.
\end{lem}

\begin{proof}
    When we consider $ \mathcal{T}_6$ modulo $5$, we get $N^2 \equiv -qe^4 \pmod{5}$. It implies 

    \begin{equation*}
         1 = \genfrac(){}{0}{-qe^4}{5}= \genfrac(){}{0}{-1}{5} \genfrac(){}{0}{q}{5} \genfrac(){}{0}{e^4}{5}= \genfrac(){}{0}{4}{5} \genfrac(){}{0}{q}{5}= \genfrac(){}{0}{q}{5}, 
    \end{equation*}
    a contradiction according to the remark \eqref{rem:3.1}. Similarly, considering the torsor $\mathcal{T}_7$, we get $\genfrac(){}{0}{q}{5}=1$, a contradiction. So neither $ \mathcal{T}_6$ nor $\mathcal{T}_7$ has any solution in $\mathbb{Z}$.
\end{proof}

Hence, combining lemmas \eqref{le:3.1}-\eqref{le:3.4}, we get $ \alpha(\Gamma)= \{1, -5pq \}$ and hence $ |\alpha(\Gamma)| = 2$. 

Now we will move to find $ |\overline{\alpha}(\overline{\Gamma})|$. As $\overline{E}_{pq}: y^2= x^3 + 20pqx$, we know that $ 1, 20pq \in \overline{\alpha}(\overline{\Gamma})$. Now all possible divisors of $b_1$ of $20pq$ modulo $ \mathbb{Q^*}^2$ is the following set

\begin{equation*}
    T=\{2,4,5,10,20, p,2p,4p,4p,5p,10p, 20p, q, 2q, 4q, 5q, 10q, 20q, pq, 2pq, 4pq, 5pq, 10pq \}.
\end{equation*}

we deliberately remove any negative values of $b_
1$ from the set $T$ as the corresponding torsors $N^2= b_1 M^4 + b_2 e^4$ will have no solutions if both $b_1$ and $b_2$ are negative. Hence, we will take all possible values of $b_1$ from the set $T$ and consider the solvability of the corresponding torsor over $\mathbb{Z}$. They are as follows:

\begin{equation*}
    \begin{split}
        \mathcal{T}_1': & N^2 = 2M^4+10pqe^4 \\
        \mathcal{T}_2': & N^2 = 4M^4+ 5pqe^4 \\
        \mathcal{T}_3': & N^2 = 5M^4+ 4pqe^4 \\
        \mathcal{T}_4': & N^2 = 20M^4+ pqe^4 \\
        \mathcal{T}_5': & N^2 = 10M^4+2pqe^4 \\
        \mathcal{T}_6': & N^2 = qM^4+20pe^4 \\
        \mathcal{T}_7': & N^2 = 2qM^4+10pe^4 \\
        \mathcal{T}_8': & N^2 = 2pM^4 + 10qe^4 \\
        \mathcal{T}_9': & N^2 = 4qM^4+5pe^4 \\
        \mathcal{T}_{10}': &  N^2 = 4pM^4+5qe^4 \\
        \mathcal{T}_{11}': & N^2 = 20pM^4+qe^4
    \end{split}
\end{equation*}

\begin{lem}
    \label{le:3.5}
    There are no integer solutions for the torsor $ \mathcal{T}_1': N^2=2M^4+10pqe^4$.
\end{lem}

\begin{proof}
    Reducing $ \mathcal{T}_1' $ modulo $5$, we get $N^2 \equiv 2M^4 \pmod 5$. Hence,

    \begin{equation*}
        1= \genfrac(){}{0}{2M^4}{5}= \genfrac(){}{0}{2}{5} \genfrac(){}{0}{M^4}{5}=  \genfrac(){}{0}{2}{5}=-1, 
    \end{equation*}
a contradiction. Thus,  $ \mathcal{T}_1' $ has no solution in $\mathbb{Z}$.
\end{proof}

\begin{lem}
       \label{le:3.6} 
      There are no integer solutions for the torsor $ \mathcal{T}_2': N^2=4M^4+ 5pqe^4$. 
\end{lem}

\begin{proof}
As $gcd(N,e)=1$ and $4$ is even, by the divisibility criteria, we can say $e$ is odd, and consequently $N$ is odd too. Thus $N^2 \equiv e^4 \equiv 1 \pmod 4$. As $ p \equiv 1 \pmod 4$ and $ q \equiv  3 \pmod 4$ by hypothesis, we have $ pq \equiv 3 \pmod 4 $. So, reducing $ \mathcal{T}_2'$ modulo $4$, we get $1 \equiv 3 \pmod 4$, a clear contradiction. Hence, $ \mathcal{T}_2'$ has no solution in $\mathbb{Z}$.

\end{proof}

\begin{lem}
    \label{le:3.7} 
     There are no integer solution for the torsors $ \mathcal{T}_3': N^2=5M^4+ 4pqe^4$ and $\mathcal{T}_4': N^2=20M^4+ pqe^4$.
\end{lem}

\begin{proof}
    Considering $ \mathcal{T}_3'$ modulo $5$, we get $N^2 \equiv 5M^4 \pmod p$. Hence,

    \begin{equation*}
        1= \genfrac(){}{0}{5M^4}{p}= \genfrac(){}{0}{5}{p}.
    \end{equation*}
    As $ p \equiv 33 \pmod {40} \equiv 3 \pmod 5$, we have $ \genfrac(){}{0}{5}{p}=-1 $. So we arrived at a contradiction and can conclude that $\mathcal{T}_3'$ has no solution in $\mathbb{Z}$. Similarly, reducing $\mathcal{T}_4'$ modulo $5$ and using Legendre symbols, we get $1= \genfrac(){}{0}{20M^4}{p} = \genfrac(){}{0}{4}{p} \genfrac(){}{0}{5}{p} \genfrac(){}{0}{M^4}{p}= \genfrac(){}{0}{5}{p}  $, a contradiction. So $\mathcal{T}_4'$ too has no solutions in $\mathbb{Z}$.
    
\end{proof}

\begin{lem}
    \label{le:3.8}
    There are no integer solutions for the torsor $ \mathcal{T}_5': N^2=10M^4+2pqe^4$.
\end{lem}

\begin{proof}
    Taking $\mathcal{T}_4'$ modulo $5$, we have $N^2 \equiv 2pq e^4 \pmod 5$. Using Legendre symbol, and remark \eqref{rem:3.1}, we get,  

    \begin{equation*}
        1= \genfrac(){}{0}{2pqe^4}{5}= \genfrac(){}{0}{2}{5} \genfrac(){}{0}{p}{5} \genfrac(){}{0}{q}{5} \genfrac(){}{0}{e^4}{5} = (-1)  (-1) (-1) 1=-1, 
    \end{equation*}
    a contradiction. So we can say that $ \mathcal{T}_5'$ has no solutions in $\mathbb{Z}$.
\end{proof}

\begin{lem}
    \label{le:3.9}
    There are no integer solutions for the torsor $ \mathcal{T}_6': N^2=qM^4+20pe^4$.
\end{lem}

\begin{proof}
    Reducing $ \mathcal{T}_6' $ modulo $5$, we get $ N^2 \equiv qM^4 \pmod {5}$. Hence,

    \begin{equation*}
        1 = \genfrac(){}{0}{qM^4}{5}= \genfrac(){}{0}{q}{5} \genfrac(){}{0}{M^4}{5} = \genfrac(){}{0}{q}{5}.
    \end{equation*}
This is a clear contradiction to the remark \eqref{rem:3.1}. Hence, it is evident that $\mathcal{T}_6'$ has no solutions in $\mathbb{Z}$.

    \end{proof}

\begin{lem}
    \label{le:3.10} 
    There are no integer solutions for the torsor $ \mathcal{T}_7': N^2=2qM^4+10pe^4$ and $\mathcal{T}_8': N^2 = 2pM^4 + 10qe^4$. 
\end{lem}

\begin{proof}
    As $N^2=2(qM^4 + 5pe^4)$, we know $N$ is even and as $gcd(N,M)=gcd(N,e)=1$, hence both $M$ and $e$ are odd. Assuming $N=2N_1$ for some $N_1 \in \mathbb{Z}$, we arrive at  $2N_1^2= qM^4 + 5pe^4$. As both $M$ and $e$ are odd, we get $M^4 \equiv e^4 \equiv 1 \pmod 8$ and using it, we can ay $ 2N_1^2 \equiv q + 5p \pmod 8$. By our assumption, $ p \equiv 33 \pmod {40} \equiv 1 \pmod 8$ and $ q \equiv 7 \pmod {40} \equiv 7 \pmod 8$. Hence, $2N_1^2 \equiv  4 \pmod 8$ which implies $N_1^2 \equiv 2 \pmod 4$, a contradiction. Thus $\mathcal{T}_7'$ has no solutions in integers. Now proceeding exactly as above, from the torsor $\mathcal{T}_8': N^2 = 2pM^4 + 10qe^4$, we get $N_1^2 \equiv 5q + p \pmod {8}$. Now using the hypothesis $ p \equiv 1 \pmod 8$ and $ q \equiv 7 \pmod 8$, we obtain  $N_1^2 \equiv 2 \pmod 4$, a contradiction. Hence, $ \mathcal{T}_8' $ too has no solutions in $ \mathbb{Z}$.
\end{proof}

\begin{lem}
    \label{le:3.11}
    There are no integer solution for the torsors $ \mathcal{T}_9': N^2=4qM^4+5pe^4$ and $\mathcal{T}_{10}': N^2=4pM^4+5qe^4$.
\end{lem}

\begin{proof}
   Considering $ \mathcal{T}_9'$ modulo $5$, we get $N^2 \equiv 4qM^4 \pmod 5$. Hence,

   \begin{equation*}
       1= \genfrac(){}{0}{4qM^4}{5} = \genfrac(){}{0}{4}{5} \genfrac(){}{0}{q}{5} \genfrac(){}{0}{M^4}{5} = \genfrac(){}{0}{q}{5}, 
   \end{equation*}
   a contradiction in accordance to the remark \eqref{rem:3.1}. So  $\mathcal{T}_8'$ has no solutions in $\mathbb{Z}$. Similarly, taking $\mathcal{T}_{10}$ modulo 5, we get $ N^2 \equiv 4pM^4 \pmod 5$, which after simplification, yields $ \genfrac(){}{0}{p}{5} =1 $. As it contradicts remark \eqref{rem:3.1}, we can say that $\mathcal{T}_{10}'$ too has no solutions in $\mathbb{Z}$.
   
\end{proof}

   \begin{lem}
       \label{le:3.12}
       There are no integer solutions for the torsors $ \mathcal{T}_{11}': N^2=20pM^4+qe^4$. 
   \end{lem}

\begin{proof}
   By reducing $\mathcal{T}_{11}'$ modulo $5$, we get $N^2 \equiv qe^4 \pmod 5$. Hence,

   \begin{equation*}
       1 = \genfrac(){}{0}{qe^4}{5} = \genfrac(){}{0}{q}{5} \genfrac(){}{0}{e^4}{5}= \genfrac(){}{0}{q}{5} , 
   \end{equation*}
   a contradiction according to remark \eqref{rem:3.1}.
\end{proof}

So from the Lemmas \eqref{le:3.5}-\eqref{le:3.12}, we can conclude that $ \overline{\alpha}(\overline{\Gamma}) = \{1, 20pq \}$ and hence $ |\overline{\alpha}(\overline{\Gamma})|=2$. 

We are ready to prove the Theorem \eqref{Thm:1.1} that is as follows.

\begin{proof}
    As $ |\alpha(\Gamma)|= 2=|\overline{\alpha}(\overline{\Gamma})| $,  from the Proposition \eqref{prop:2.1}, we have

    \begin{equation}
        2^{r}= \frac{1}{4}( |\alpha(\Gamma)| \cdot |\overline{\alpha}(\overline{\Gamma})|=1,
    \end{equation}
    from which we get that $r=r(E_{pq})=0$.
\end{proof}

The following table \eqref{tab:1} confirms the result of Theorem \eqref{Thm:1.1}. Here we take some pairs of primes $p$ and $q$ which satisfy the condition of the Theorem \eqref{Thm:1.1} and list down the rank of the corresponding elliptic curve $E_{pq}$. All computations are done using SAGE \cite{SA}.

\begin{table}[ht]
    \centering

\begin{tabular}{|c| c| c|} 

\hline
 $p$ & $q$ & $\text{rank of}~ E_{pq}$ \\
\hline 
73 & 7 & 0 \\
\hline
113 &  7 & 0 \\
\hline
73 & 47 & 0 \\
\hline 
113 & 47 & 0 \\
\hline
73 & 127 & 0 \\
\hline 
113 & 127 & 0 \\ 
\hline

\end{tabular}
 \caption{The values of $p$ and $q$ that satisfies the condition of Theorem \eqref{Thm:1.1}}
    \label{tab:1}
\end{table}

\begin{cor}
    \label{cor:3.1} Under the same assumption as of Theorem \eqref{Thm:1.1}, if the congruence $q \equiv 7 \pmod {40}$, is replaced by $ q  \equiv 7 \pmod {20}$, then the rank of elliptic curve $E_{pq}$ given by equation \eqref{eq:1.1} is at most one. 
\end{cor}

\begin{proof}
Throughout all lemmas \eqref{le:3.1}-\eqref{le:3.12}, we have used the condition that $ q \equiv 2 \pmod{5}$, but not  $ q \equiv 7 \pmod 8$ except lemma \eqref{le:3.10}. That particular lemma uses the condition $ q \equiv 7 \pmod 8$ and hence the torsors $\mathbf{T}_7'$ and $ \mathcal{T}_8'$ did use the same assumption. According to this new assumption, $q$ does not satisfy $ q \equiv 7 \pmod {8}$.  This implies torsors $\mathbf{T}_7'$ and $ \mathcal{T}_8'$ may have integer solutions under this new condition. So from lemma \eqref{le:3.10}, we can say that it is possible $ 2p, 10p, 2q, ~\text{and}~ 10q \in \overline{\alpha}(\overline{\Gamma})$, making $ |\overline{\alpha}(\overline{\Gamma})| \leq 6$. As a consequence, we obtain $ o \leq r \leq 1$
\end{proof}

Using this corollary, we will now prove the Theorem \eqref{thm:1.2}

\begin{proof}
If $ p = 40 k + 33$ and $q = 40l + 27$ for some $k,l \in \mathbb{Z}$, then it satisfies the assumptions of the corollary \eqref{cor:3.1}. So the torsors that may have integer solutions are $ \mathcal{T}_{7}'$ and $\mathcal{T}_{8}'$. Now 

\begin{equation*}
    \mathcal{T}_{8}': N^2 = 2(33 + 40k) + 10 ( 40l +27)= 16(5k+25l +21l),
\end{equation*}

implies $ \mathcal{T}_{8}'$ has a solution $(N,M,e)= ( 4 \sqrt{5k+25l +21}, 1,1)$. Hence, $ 2p, ~\text{and}~ 10 q \in |\overline{\alpha}(\overline{\Gamma})| $. Now if the torsor $ \mathcal{T}_{7}'$ has an integer solution, it would imply $ 2q, 10p \in  |\overline{\alpha}(\overline{\Gamma})|$ which in turn gives $  |\overline{\alpha}(\overline{\Gamma})| = 6$. Now from the Proposition \eqref{prop:2.1}, we obtain

\begin{equation*}
    2^{r}= \frac{1}{4} |\alpha(\Gamma)| |\overline{\alpha}(\overline{\Gamma})|= 3,
\end{equation*}

which is impossible. Hence,  $ |\overline{\alpha}(\overline{\Gamma})|=4$ and as a consequence, we obtain $ r(E_{pq})=1$.

\end{proof}

Here also, to confirm our hypothesis in Theorem \eqref{thm:1.2}, we list down some pairs of $p$ and $q$ that satisfy the conditions of Theorem \eqref{thm:1.2} and their corresponding rank in Table \eqref{tab:2}.

\begin{table}[ht]
    \centering

\begin{tabular}{|c| c| c|c| c|} 

\hline
$k$ & $l$ & $p$ & $q$ &  $\text{rank of}~ E_{pq}$ \\
\hline 
2 & 2 & 113 & 107 & 1 \\
\hline
2 &  5 & 113 & 227 & 1\\
\hline
2 & 8 & 113 & 347 & 1 \\
\hline 
7 & 5 & 313 & 227 & 1 \\
\hline
10 & 5 & 433 & 227 & 1 \\
\hline 
7 & 8 & 313 & 347 & 1 \\ 
\hline

\end{tabular}
 \caption{The values of $p$ and $q$ that satisfies the condition of Theorem \eqref{thm:1.2}}
    \label{tab:2}
\end{table}


Now using Theorem \eqref{Thm:1.1}, \eqref{thm:1.2} and  \eqref{thm:2.2}, we get the following corollary.

\begin{cor}
 \label{cor:3.2}
 Consider the field $\mathbb{K}=\mathbb{Q}(i)$. Then we can say the following about the rank of the elliptic curve $E_{pq}$ given by \eqref{eq:1.1} over $\mathbb{K}$. \\

 (i) If $ p$ and $q $ are distinct odd primes satisfying the condition as in Theorem \eqref{Thm:1.1}, then the rank of $E_{pq} $ over $K$ is zero. \\

 (ii) If $ p$ and $q $ are distinct odd primes, and $l,k \in \mathbb{Z}$  such that they satisfy the hypothesis of Theorem \eqref{thm:1.2}, then the rank of $E_{pq} $ over $\mathbb{K}$ is $2$.

\end{cor}

\begin{proof}
    (i) If $\mathbb{K}=\mathbb{Q}(i)$, then from Theorem \eqref{thm:2.2}, we have $m=-1$, and, hence $E_{pq}[-1]: y^2 =x^3 -(-1)^2 ~5pqx= x^3 -5pqx= E_{pq}$. Hence, using the Theorem \eqref{thm:2.2}, we can say 

    \begin{equation*}
        \begin{split}
            rank(\mathbb{K}) &= rank(E_{pq}(\mathbb{Q})) + rank(E_{pq}[-1](\mathbb{Q}))\\
            & = 2~ rank(E_{pq}(\mathbb{Q})).
        \end{split}
    \end{equation*}
From Theorem \eqref{Thm:1.1}, we know that $rank(E_{pq}(\mathbb{Q}))=0$. So our claim follows directly by using the above equation.

(ii) From Theorem \eqref{thm:1.2}, we know that $ rank(E_{pq}(\mathbb{Q}))=1$. So  by using the above equation $   rank(\mathbb{K}) = 2~ rank(E_{pq}(\mathbb{Q}))$, we get the desired conclusion. 
\end{proof}

\section{Torsion of $\texorpdfstring{E_{pq}}{}$}

From the Theorem \eqref{thm:2.1}, we know that either $Y=0$ or $Y^2 \mid 500p^3q^3$. So when $Y=0$, then from the equation \eqref{eq:1.1}, we have either $X=0$ or $X^2 = 5pq$. Now when $X=0$, we get the Torsion point $P=(0,0)$ of order $2$. Now if $X^2 = 5pq$, then clearly $X$ is odd and hence $X^2 \equiv 1 \pmod 8$. However, from hypothesis, $p \equiv 33 \pmod {40} \equiv 1 \pmod 8$ and $ q \equiv 7 \pmod 8$, and hence $5pq \equiv 3 \pmod 8$ which is a contradiction. So $X^2 =5pq $ has no integral solution. 

\subsection{Proof of the Theorem $\texorpdfstring{\eqref{thm:1.3}}{}$}

\begin{proof}
    Let us assume $P=(X, Y)$ be the torsion point of the elliptic curve given by equation \eqref{eq:1.1}. From the above discussion, we know $P=(0.0)$ is a point of order $2$ on the given elliptic curve. We claim that $E_{pq}$ has no other non-trivial torsion point except $P =(0,0)$ which implies that $Y^2 \mid 500p^3q^3$ gives no other torsion point. We prove this in the following way:

    Suppose our assumption was wrong. Then $ Y^2 \mid 500p^3q^3 $ implies 
    \begin{equation*}
         Y^2 \in S= \{ 1, 4, 25, 100, p^2, (2p)^2, (5p)^2, (10p)^2,  q^2, (2q)^2, (5q)^2, (10q)^2, (pq)^2, (2pq)^2, (5pq)^2, (10pq)^2 \}. 
    \end{equation*}
    Now we will consider each different scenario.

    \textit{Case-I}: Let $Y^2=1$. Then from the equation \eqref{eq:1.1}, we have $1=X(X^2-5pq)$. Hence $X=(X^2-5pq)= \pm 1$. Now if $X=1$, then $X^2-5pq=1$ which in turn gives $5pq=0$, a contradiction. Now when $X=-1$, we get $5pq=2$, a contradiction.

    \textit{Case-II}: Now we assume $Y^2=4$. Hence, from the equation \eqref{eq:1.1}, we get $X(X^2-5pq)=4$. Now we will consider a couple of different scenarios:

    \textit{Sub-case I}: Let $X=\pm 1$. Then we have $X^2 -5pq = \pm 4$. So if $X=1$, then we have $ 5pq =-3$, a contradiction. Similarly, when $X=-1$, we have $ pq=1$, contradicting our hypothesis.

    \textit{Sub-case II}: Let $X = \pm 2$ which in turn implies $X^2 - 5pq = \pm 2$. After simplification, we have $ 5pq= -2 ~\text{or} -6$ whenever $ X = 2$ or $X =-2$ respectively. So, we get $ 5 \nmid 2$ or $5 \nmid 6 $, a contradiction.

    \textit{Sub-case III}: Assume $X = \pm 4$.  Then $X^2 -5pq = \pm 1$, which in turn implies $ pq =-3$ when $X = 4$ and $5pq= -17$ when $X=-4$. From the first instance, it is clear that one of the $p$ or $q$ has to be $\pm 1$ which is a contradiction. The second scenario is impossible too as $ 5 \nmid 17$.

    \textit{Case-III}: Let us now focus on the scenario when $Y^2 = 25$. Then, we have $25= X(X^2-5pq)$. Now we split this case into three different scenarios.

    \textit{Sub-case I}: Let us take $X= \pm 1$. Then we have $5pq = -24$ or $5pq = 26$, a contradiction to the hypothesis as both $p$ and $q$ are odd, and hence $5pq $ is odd too.

     \textit{Sub-case II}: When $X = \pm 5$, we got $pq =-4$ or $pq=6$. It implies that at least one of $p$ and $q$ is even, a contradiction to our hypothesis.

     \textit{Sub-Case III}: Now if we assume $X= \pm 25$, we get $ 5pq = 624$ or $ -626$, a clear contradiction as $5pq \equiv 1 \pmod 2$  odd and both $624 $ and $ 626$ are equivalent to $ 0 \pmod 2$.

     \textit{Case-IV} Now we will turn our attention to the scenario where $Y^2 = p^2 q^2$. Hence, from the equation \eqref{eq:1.1}, we have $ (pq)^2 = X (X^2-5pq)$. We now will divide this into several sub-cases as follows.

     \textit{Sub-case I}: Let $X= \pm 1$ which implies $ X^2 -5pq = \pm 1$. After simplification, we have $ 5pq = 1 -(pq)^2$ or $ 5pq= 1+ (pq)^2$. As $ 5pq \equiv 3 \pmod 4$ and $ (pq)^2 \equiv 1 \pmod 4$, we have $ 1 + (pq)^2 \equiv 2 \pmod 4 $ and $ 1- (pq)^2 \equiv 0 \pmod 4$. Hence, none of the equations  $ 5pq = 1 -(pq)^2$ and  $ 5pq= 1+ (pq)^2$ has any integral solution.

     \textit{Sub-Case II}:
     $ X= \pm (pq)^2$. Then $X^2 - 5pq = \pm 1 $  which can be simplified to $ pq( p^3 q ^3 -5) = \pm 1$. It implies that $pq$ has to be either $1$ or $-1$ in any case. This is a contradiction to the hypothesis that both $p$ and $q$ are odd primes.

      \textit{Sub-Case III}:  $X = \pm pq$, which gives $ X^2 -5pq = \pm pq$. After implications, we have $ pq =6$ or $ pq =4$, which implies that one of $p$ or $q$ has to be even, a contradiction.

      \textit{Sub-Case IV}: $X = \pm p$. Then $ X^2 -5pq = \pm pq^2$, which reduces to $ p = q( 5 \pm q)$. As $ q \equiv 1 \pmod 8$, $ 5 \pm q \equiv 4,6 \pmod 8 $ and hence $ q( 5 \pm q) \equiv 4, 2 \pmod 8$. As $ p \equiv 1 \pmod 8$, we arrived at a condition here too. The case of $X = \pm  q$ can be discarded similarly. 

      \textit{Sub-Case V}: $X = \pm p^2$. Now $ X^2 -5pq = \pm q^2 $, which after simplification, yields $ p^4 = q(5p \pm q)$.  As $ (5p \pm q) \equiv 4, 2 \pmod 8$ and $ p^4 \equiv 1 \pmod 8$, we get a contradiction. The case $ X = \pm q^2 $ can be discarded similarly.

      \textit{Sub-Case VI}: If $X = \pm p^2 q$, then we get $ p^4 q = 5p \pm 1$. As the right-hand side is equivalent to $ 7 \pmod 8$ whereas the left-hand side is equivalent to $ 4, 6 \pmod 8$, we get a contradiction.  

      Similarly, as above, we can show that none of the elements of the set $S$ can be a possible value of $Y^2$ and consequently the value of $X$. Hence, torsion subgroup of $E_{pq} \simeq \mathbb{Z}/2 \mathbb{Z}$.     
\end{proof}

\section{Concluding Remarks}

In this article, we have found conditions on the distinct primes $p$ and $q$ such that the rank of the elliptic curve $E_{pq}$ given by the equation \eqref{eq:1.1} over $\mathbb{Q}$ is either zero or one. However, can we find conditions on $p$ and $q$ such that the rank of the elliptic curve given by \eqref{eq:1.1} is two or higher? It is an interesting question and we encourage readers to pursue this question.

\section*{Funding and Conflict of Interests/Competing Interests} The author has no financial or non-financial interests to disclose that are directly or indirectly related to the work. The author has no funding sources to report.

\section*{Data availability statement} No outside data was used to prepare this manuscript.

\end{document}